\documentclass[a4paper]{article}
\usepackage[comma,numbers,square,sort&compress]{natbib}
\usepackage{amsmath,amssymb,bm,amsthm}
\theoremstyle{definition} 
\usepackage{subfigure,pgf,tikz,pst-node}
\usetikzlibrary{arrows,shapes,automata,positioning,fit}
\renewcommand{\bar}{\overline}
\renewcommand{\top}{\intercal}

\newcommand{\e}{\varepsilon}
\newcommand{\R}{\mathbb{R}}

\usepackage{bm} 
\newcommand{\N}{\mathcal{N}}

\newtheorem{theorem}{Theorem}
\newtheorem{lemma}[theorem]{Lemma}

\newtheorem{remark}[theorem]{Remark}

\newtheorem{ass}{Assumption}

\newcommand{\pb}{\noindent\textbf{Proof. } }
\newcommand{\pe}{\hfill\rule{4pt}{8pt}}

\usepackage{pgf,tikz}
\usetikzlibrary{arrows,shapes,automata,positioning,fit}
\usepackage{pst-node}

\begin{document}
	
	\title{Distributed Optimization with Inexact Oracle}
	
	\author{Kui Zhu, Yichen Zhang, and Yutao Tang \footnote{This work was supported by National Natural Science Foundation of China under Grants 61973043. K. Zhu, Y. Zhang, and Y. Tang are with the School of Artificial Intelligence, Beijing University of Posts and Telecommunications, Beijing 100876. P.\,R. China. Emails: 1239366542@qq.com, zhangyc930@163.com, yttang@bupt.edu.cn}}
	
	\date{}
	
	\maketitle
	
	{\noindent\bf Abstract}: In this paper, we study the distributed optimization problem using approximate first-order information. We suppose the agent can repeatedly call an inexact first-order oracle of each individual objective function and exchange information with its time-varying neighbors. We revisit the distributed subgradient method in this circumstance and show its suboptimality under square summable but not summable step sizes.  We also present several conditions on the inexactness  of the local oracles to ensure an exact convergence of the iterative sequences towards the global optimal solution.  A numerical example is given to verify the efficiency of our algorithm.
	

\section{Introduction}

Consider a network of $m$ agents and each agent $i$ corresponds to a convex objective function $f_i(x)$. Suppose the agent can collect the data of each individual function and exchange some information with other agents. We aim at an optimal solution and the optimal value of the optimization problem:
\begin{align}\label{prop:main}
\min_{x\in \R^n }f(x)=\sum_{i=1}^m f_i(x)
\end{align}
This problem has been coined as distributed optimization and shown to have many applications in multi-robot systems, smart grid, communication networks, and large-scale machine learning \cite{nedic2018distributed, yang2017distributed}.  

Among different classes of numerical algorithms for this problem, methods that use the first-order information (i.e., the function value $f_i(s)$ and its subdifferential $\partial f_i(s)$ at specific point $s\in \R^n$) have been intensively studied due to the cheap iteration cost and well-established convergence properties. In the seminal work \cite{nedic2009distributed}, a subgradient rule was proposed to solve the distributed optimization problem, where each agent performs a consensus step and then a descending step along the subgradient direction of each local objective function. Since then, many attempts have been made to solve such distributed optimization problems in different settings by primal domain, dual domain, or primal-dual methods, just to list a few \cite{duchi2011dual, shi2013reaching, yi2015distributed, shi2015EXTRA, xi2016distributed, lei2016primal, zeng2017distributed}. 
Note that these optimization results can be roughly regarded as distributed extensions of the classical centralized first-order methods and all require the exact first-order information of the individual objective functions. 

In practice, the first-order information of the objective functions (especially the subgradients) may not be exactly available without knowing the function structure. It is usually calculated by considering another auxiliary optimization problem, which is often solved with certain accuracy \cite{correa1993convergence, bertsekas2015convex}. This issue will be aggravated by the possible environment noise \cite{nedic2010effect}. Thus it is crucial to consider the solvability of the given optimization problem using only approximate first-order information. In fact, this issue has been paid much attention to in the field of mathematical programming. For example, the well-known notion of $\varepsilon$-subgradient has been delivered and intensively used when the exact subgradient is not available \cite{polyak1987introduction, nedic2010effect, bertsekas2015convex, rasch2020inexact}. 

Recently, the authors in  \cite{devolder2014first} introduced a very general notion of inexact first-order oracle to describe the inexactness of the first-order information for given objective functions. It has been intensively utilized to analyze the behavior of several typical centralized inexact first-order methods in the past few years and also been applied to study distributed optimization algorithms \cite{jakovetic2014fast, nedic2017achieving, qu2018harnessing}. Typically, these papers are often focused on the evolution of some averaging dynamics while the oracle inexactness is interpreted from the coordinate error among agents to prove the convergence performance of the proposed distributed algorithms therein. Note that they still assume the exact first-order information of each individual objective function. 
With the aforementioned observation, a natural question arises: {\it Whether and how can we develop distributed first-order optimization methods with local inexact oracles to solve the optimization problem \eqref{prop:main}?}

Having this question in mind, we will focus the optimization problem \eqref{prop:main} and seek effective distributed algorithms to solve it providing that only approximate first-order information of the objective functions is available. To this end, we revisit the well-known distributed subgradient method proposed in \cite{nedic2009distributed} but working with an inexact oracle of each $f_i$. We first show the convergence of this rule under diminishing step sizes and then investigate its optimality or suboptimality under different oracle conditions. 

The main contribution can be summarized as follows. 

\begin{itemize}
	\item We present a novel distributed algorithm to solve the consensus optimization problem under time-varying communication graphs. This rule can be taken a distributed counterpart of existing primal gradient method with inexact oracle in \cite{devolder2014first}. Compared with existing distributed algorithms in \cite{nedic2009distributed, duchi2011dual, jakovetic2014linear, shi2015EXTRA, xi2016distributed, lei2016primal}, this rule can work with only approximate first-order information of the objective functions. 
	\item We establish the convergence of our proposed algorithm under diminishing step sizes and show that it can approximately solve the problem provided that each local objective function has an inexact first-order oracle. Moreover, we specify a group of conditions under which the exact convergence of such iterative sequences towards the global optimal solution can be ensured even with inexact local oracles. 
\end{itemize}

The rest of this paper is organized as follows. In Section \ref{sec:pre}, we first introduce  some notations and concepts for the following analysis. Then, we present a novel distributed rule with inexact oracle in Section \ref{sec:algorithm}. After that, we consider the convergence and suboptimality issues of our algorithm under different conditions in Section \ref{sec:convergence}. We also give a simulation example in Section \ref{sec:simu} to  illustrate the effectiveness of our design. Finally, we concludes the paper with some remarks in Section \ref{sec:con}.

\section{Preliminary}\label{sec:pre}

We briefly present some preliminaries about our notations and terminology for the following analysis \cite{devolder2014first, bertsekas2015convex}. 

\subsection{Notation}
Let $\R^n$ be the $n$-dimensional Euclidean space and $\R^{n\times m}$ be the set of all $n\times m$ matrices with all entries in $\R$.  A vector is viewed as a column, unless clearly stated otherwise.  ${\bf 1}_n$ (or ${\bf 0}_n$) denotes an $n$-dimensional all-one (or all-zero) column vector and ${\bm 1}_{n\times m}$ (or ${\bm 0}_{n\times m}$) all-one (or all-zero) matrix. 
$\mbox{diag}(b_1,\,{\dots},\,b_n)$ denotes an $n\times n$ diagonal matrix with diagonal elements $b_1,\,\dots,\, b_n$.  $\mbox{col}(a_1,\,{\dots},\,a_n) = {[a_1^\top,\,{\dots},\,a_n^\top]}^\top$ for column vectors $a_i\; (i=1,\,{\dots},\,n)$.  For a vector $x$ (or matrix $A$) , $||x||$ ($||A||$) denotes its Euclidean (or spectral) norm. When $A$ is square,  its trace and Frobenius norm are denoted by $\mbox{Tr}(A)$ and $||A||_{\rm F}=\mbox{Tr}(A^\top A)$, respectively. 

A column (or row) vector $x$ is said to be a stochastic vector when its components are nonnegative and their sum is equal to $1$. A square matrix is doubly stochastic when its rows are stochastic vectors, and its columns are also stochastic vectors.

\subsection{Convex function and inexact oracle}

Consider a function $f\colon \R^m  \to \R $. It is said to be convex if for  $0\leq a \leq 1$ and  $\forall \zeta_1,\zeta_2 \in \R^m$, we have
\begin{align*}
	f(a\zeta_1+(1-a)\zeta_2)\leq a f(\zeta_1)+(1-a)f(\zeta_2)
\end{align*}
If it is differentiable, the convexity of this function is equivalent with the following inequality: 
$$f(\zeta_1)\geq f(\zeta_2)+\nabla f(\zeta_2)^T(\zeta_1-\zeta_2),\, \forall \zeta_1,\zeta_2 \in \R^m$$

Function $g\colon \R^m\to \R^m$ is Lipschitz with constant $\vartheta_1>0$ if 
\begin{align*}
	||g(\zeta_1)-g(\zeta_2)||\leq \vartheta_1 ||\zeta_1-\zeta_2||,\, \forall \zeta_1,\, \zeta_2 \in \R^m
\end{align*}

A convex function $f$ is said to admit a first-order $(\delta,\, L)$-oracle if for any given point $x\in \R^m$, we have a pair $(f_{\delta,\,L}(x),\, g_{\delta,\,L}(x))\in \R\times \R^m$ such that
\begin{align*}
	0\leq f(y)-f_{\delta,\,L}(x)-g_{\delta,\,L}(x)^\top (y-x)\leq \frac{L}{2}\|y-x\|^2+\delta
\end{align*}
holds for any $y\in \R^m$ with $\delta \geq 0$ being the accuracy of the inexact first-order oracle.  If $f$ admits a $(\delta,\, L)$-oracle, then $cf$ admits a  $(c\delta,\,cL)$-oracle for any constant $c>0$. If $f_i$ admits a $(\delta_i,\, L_i)$-oracle, $i=1,\,\dots,\,m$, then $\sum_{i=1}^m f_i$ admits a $(\sum_{i=1}^m\delta_i,\,\sum_{i=1}^m L_i)$-oracle.

When a differentiable convex function $f$ has  $L$-Lipschitz gradients, it admits a $(0,\,L)$-oracle with $f_{0,\,L}(x)=f(x)$ and $g_{0,\,L}(x)=\nabla f (x)$. This is the exact first-order information of function $f$. When a (nonsmooth) convex function $f$ can be well approximated by a smooth convex function $\bar f$ with Lipschitz gradients in the sense that their difference is bounded, the exact first-order information of function $\bar f$ provides an inexact oracle for function  $f$.

\section{Distributed Rule with Inexact Oracle}\label{sec:algorithm}

The goal of this multi-agent system is to collectively minimize the global objective function $f$ in a distributed manner. As usual, we assume that the optimal value $f^*=\min_{x\in \R^n} f(x)$ is finite and the optimal solution set $\mathcal{X}^*=\{x \in \R^n \mid f(x)=f^*\}$ is nonempty.  
 
In this paper, we are more interested in the solvability of the distributed optimization problem working with only an inexact oracle of the local objective function. 

\begin{ass}\label{ass:inexact}
	For each $i=1,\,\dots,\,m$, function $f_i$ is convex and admits a first-order $(\delta,\,L)$-oracle.
\end{ass}

Several nontrivial classes of examples satisfying this assumption can be found in \cite{devolder2014first}. Although smooth convex functions with Lipschitz gradients naturally satisfy this assumption, $f_i$ needs not to be differentiable at all. We will give a nonsmooth example later in Section \ref{sec:simu}. 
 
We suppose agent $i$ maintains an estimate $x_i(k)$ of the global optimal solution to problem \eqref{prop:main} at time instant $k=1,\,2,\,\dots$ and can share its estimate with other agents. Meanwhile, it will collect the data of problem \eqref{prop:main} by calling the first-order $(\delta,\,L)$-oracle of $f_i$ at time $k$. Suppose the  oracle will return a pair of $(f_{i,\,\delta}(k),\, g_{i,\,\delta}(k))$ satisfying the following assumption:  
\begin{ass}\label{ass:boundedness}
	There exists a constant $C>0$ such that $\max_i\{\|g_{i,\delta}(k)\|\}\leq C$ for all $k\geq0$ and $i=1,\,\dots,\,m$.
\end{ass}

Each agent will update its estimate based on the inexact first-order information of the local objective function and all available estimates including its own and those received from others.  To describe the communication topology, we introduce a sequence of undirected graphs $\{\mathcal{G}_k\}$ with node set $\N=\{1,\,\dots,\,m\}$ and edge set $\mathcal{E}_k$. An edge $(i,\,j)\in \mathcal{E}_k$ means agent $i$ and agent $j$ can communicate with each other at time instant $k$. Motivated by the well-known distributed subgradient method in \cite{nedic2009distributed}, we consider the following iteration with inexact oracles of the individual functions:  
\begin{align}\label{alg:main}
	x_i(k+1)&=\sum_{j=1}^m w_{ij}(k)x_j(k)-\alpha_k g_{i,\delta}(k)
\end{align}
where $g_{i,\,\delta}(k)$ is an inexact subgradient of function $f_i$ at point $x_i(k)\in \R$ by calling its inexact oracle and $\alpha_k>0$ is the step size to be specified later. Here $w_{ij}(k)$ is a nonnegative weight with the property that $w_{ij}(k)>0$ if and only if node $i$ receives information from node $j$ at time $k$. We use a matrix $W(k)=[w_{ij}(k)]_{m\times m}$ to represent these weights and denote  ${\bf x}(k)=\mathrm{col}(x_1(k),\,\dots,\,x_m(k))$, $g_{\delta}({\bf x}(k))=\mathrm{col}(g_{1,\delta}(k),\,\dots,\,g_{m,\delta}(k))$ for short. Then we can rewrite the distributed rule \eqref{alg:main} into a compact form as follows:
\begin{align}\label{alg:main-compact}
	{\bf x}(k+1)=[W(k)\otimes I_n]{\bf x}(k)-\alpha_k g_{\delta}({\bf x}(k))
\end{align}
We may denote $g_{\delta}(k)=g_{\delta}({\bf x}(k))$ for short.  This rule \eqref{alg:main} extends the distributed subgradient method developed in \cite{nedic2009distributed} to the case when only approximate first-order information of the local objective function is available. 

Here are some assumptions on the communication graphs to guarantee a consensus when $g_{\delta}({\bf x}(k))={\bm 0}$, which have been widely used in literature \cite{nedic2009distributed,nedic2010constrained,liu2017convergence}.
\begin{ass}\label{ass:graph} 
	Graph $\mathcal{G}_k$ is connected for any $k\geq 1$.
\end{ass}
\begin{ass}\label{ass:graph-weight}
	For any $k\geq 1$,  the matrix $W(k)$ is doubly stochastic. Additionally, there is a constant $\eta>0$ such that, for any $i=1,\, \dots,\, m$,  $w_{ii}(k)\geq \eta$ and $w_{ij}(k) \geq \eta$ if $(j,\,i) \in \mathcal{E}_k$.
\end{ass}

In the following section, we will study the convergence of our algorithm \eqref{alg:main} and  show the (approximate) solvability of distributed optimization problem \eqref{prop:main} by this rule with the inexact oracle.  For the ease of presentation, we assume $n=1$ and multiple dimensional analysis can be done for each component without any technical obstacles. 

\section{Convergence Analysis}\label{sec:convergence}

This section is dedicated to the performance analysis of the distributed rule \eqref{alg:main}.  Denote the average of all agents' estimates at time $k$ by $x_{\rm av}(k)=\frac{{\bf 1}^\top{\bf x}(k)}{m}$ and let $\bar {\bf x}(k)={\bf x}(k)-{\bf 1}x_{\rm av}(k)$. It follows then
\begin{align}\label{sys:error}
	\bar {\bf x}(k+1)&=Q(k) \bar {\bf x}(k)-\alpha_k R g_{\delta }(k)
\end{align}
where $Q(k)=W(k)-\frac{{\bf 1}{\bf 1}^\top}{m}$ and $R=I_m-\frac{{\bf 1}{\bf 1}^\top}{m}$.

We denote  $\Phi(k,\,s)\triangleq \prod_{k'=s}^k W(k')$.  The following lemma ensures an asymptotic consensus of our algorithm \eqref{alg:main} under diminishing step sizes.

\begin{lemma}\label{lem:consensus}
	Suppose Assumptions \ref{ass:inexact}--\ref{ass:graph-weight} hold and choose the step size $\{\alpha_k\}$ such that  $\lim_{k\to \infty} \alpha_k=0$. Along the trajectory of algorithm \eqref{alg:main}, all estimates of agents can achieve a consensus on their time-varying average, that is, $\bar {\bf x}(k) \to {\bf 0}$ as $k\to \infty$. 
\end{lemma}
\pb
To prove this lemma, we view it as a perturbed stability problem of \eqref{sys:error} with a perturbation $\alpha_k Rg_{\delta}(k) $.   Note that $\prod_{k'=s}^k Q(k')=\Phi(k,\,s)-\frac{{\bf 1}{\bf 1}^\top}{m}$. It follows then
\begin{align*}
	\bar {\bf x}(k+1)&= \prod_{k'= K}^k Q(k') \bar {\bf x}(K)-\sum_{k'=K}^{k-1} \alpha_{k'} \prod_{i=k'+1}^{k} Q(i) Rg_{\delta }(k') - \alpha_k  Rg_{\delta }(k)\\
	&=[\Phi(k,\,K)-\frac{{\bf 1}{\bf 1}^\top}{m}] \bar {\bf x}(K) - \alpha_k  Rg_{\delta }(k) -\sum_{k'=K}^{k-1} \alpha_{k'} [\Phi(k,\, k'+1)-\frac{{\bf 1}{\bf 1}^\top}{m}] Rg_{\delta }(k')
\end{align*}
for any $k>K$. We will show that for any $\e>0$, there exist a large enough $K$ such that $ ||\bar {\bf x}(k+1)||<\e $ holds for any $k>K$.

Under the theorem conditions, we first recall Lemma 1 in \cite{nedic2018distributed}. It follows that
	\begin{align*}
		([\Phi(k,\,s)]_{ij}-\frac{1}{m})^2 \leq(1-\frac{\eta}{2m^2})^{k-s}
	\end{align*}
holds for any $s>0$ and $k\geq s$.  This jointly with the fact $\sqrt{1-\mu}\leq 1-\frac{\mu}{2}$ when $0<\mu<1$ implies that 
	\begin{align*}
		|[\Phi(k,\,s)]_{ij}-\frac{1}{m}|&\leq\sqrt{(1-\frac{\eta}{2m^2})^{k-s}} \leq (1-\frac{\eta}{4m^2})^{k-s}
	\end{align*} 
Denote $q=1-\frac{\eta}{4m^2}$ for short. Since $||[\Phi(k,\, s)-\frac{{\bf 1}{\bf 1}^\top}{m}]||_{\rm F}\leq m \max_{i,\,j}\{ |[\Phi(k,\,s)]_{ij}-\frac{1}{m}|\}\leq mq^{k-s}$.  Noticing that $||g_{\delta }(k)||\leq \sqrt{m}C$ under Assumption \ref{ass:boundedness}, one can obtain   
	\begin{align*}
		||\bar {\bf x}(k+1)||&\leq m q^{k-K} ||\bar {\bf x}(K)||+ \alpha_k \sqrt{m} C||R||+m^{3/2}  C||R||\sum_{k'=K}^{k-1} \alpha_{k'} q^{k-k'-1}\\
		&\leq  m q^{k-K} ||\bar {\bf x}(K)||+ \sup_{k'\geq K} \{\alpha_{k'}\}  \frac{m^{3/2} C||R||}{q(1-q)}
	\end{align*}

Since $\lim_{k\to \infty} \alpha_k=0$, for a given $\e>0$, there exists a large enough $K'$ such that  $\sup_{k'\geq K'} \{\alpha_{k'}\}\leq  \frac{q(1-q)\e }{2m^{3/2}C||R||}$. Note that there exists an integer $K>K'$ such that $q^{k-K'}<\frac{\e}{2m||\bar {\bf x}(K')||}$ when $k>K$ due to the fact that $0<q<1$. Therefore, we have $||\bar {\bf x}(k+1)||<\e$ for any $k>K$. The proof is thus complete.
\pe

Lemma \ref{lem:consensus} states that all agents' estimates of the global optimal solution eventually reach a consensus if $\lim_{k\to \infty}\alpha_k=0$. To establish the convergence of \eqref{alg:main}, we shall focus on the evolution of   $x_{\rm av}$ governed by
\begin{align}\label{sys:average}
	\begin{split}
		x_{\rm av}(k+1)&=x_{\rm av}(k)-\alpha_k \frac{{\bf 1}^\top g_{\delta}(k)}{m}
	\end{split}
\end{align}

Denote $\Delta_i(k)=g_{i,\delta}(k)(x_i(k)-x_{\rm av}(k))$. By the definition of $(\delta,\,L)$-oracle, we have
\begin{align*}
	& f_i(y_i)-f_{i,\delta}(k)-g_{i,\delta}(k) (y_i-x_{\rm av}(k))+\Delta_i(k)\geq 0 \\
	&  f_i(y_i)-f_{i,\delta}(k)-g_{i,\delta}(k)(y_i-x_{\rm av}(k))+\Delta_i(k)\\
	&\qquad \qquad \quad   \leq L\|y_i-x_{\rm av}(k)\|^2+L\|x_{\rm av}(k)-x_i(k)\|^2+\delta
\end{align*}
Summing up these inequalities together, we have
\begin{align}\label{eq:inexact-distributed}
\begin{split}
\Xi ({\bf y}, k)&\geq 0 \\
\Xi ({\bf y}, k)&\leq L\|{\bf y}-{\bf 1}  x_{\rm av}(k)\|^2+L\| {\bar {\bf x}}(k)\|^2+m\delta
\end{split}
\end{align}
where ${\bf y}=\mathrm{col}(y_1\,\dots,\,y_m)$, $\Xi ({\bf y}, k)\triangleq \sum_{i=1}^mf_i(y_i)-(\sum_{i=1}^mf_{i,\delta}(k)-\sum_{i=1}^m \Delta_i(k))- \sum_{i=1}^m g_{i,\delta}(k)^\top (y_i-x_{\rm av}(k))$. 

It is interesting to remark that the above two inequalities can be interpreted in the terminology of inexact oracle. In fact, by taking ${\bf y}={\bf 1}y$  for any $y \in \R$, we can find that $ \sum_{i=1}^mf_{i,\delta}(k)- \sum_{i=1}^m \Delta_i(k)$ and $ \sum_{i=1}^m g_{i,\delta}(k)$ are the approximate first-order information of function $f$ at the point $x_{\rm av}(k)$ with time-varying oracle accuracy $L\|{\bar {\bf x}(k)}\|^2+m\delta$.  Since $\|{\bar {\bf x}(k)}\|$ will eventually converge to zero under a diminishing step size, we are going to get an inexact oracle of $f$ with accuracy $m \delta$ in an asymptotic sense. This fact has been extensively exploited in \cite{jakovetic2014fast,nedic2017achieving,qu2018harnessing} to establish the convergence of the distributed optimization algorithms therein when $\delta=0$. 

With these observations, it is time to present the first main theorem of the paper. 

\begin{theorem}\label{thm:inexact}
Suppose Assumptions \ref{ass:inexact}--\ref{ass:graph-weight} hold and choose a positive step size $\{\alpha_k\}$ to satisfy the condition
	\begin{align*}
		\sum_{k=1}^\infty \alpha_k=\infty,\quad \sum_{k=1}^\infty \alpha_k^2<\infty
	\end{align*}
Then, the following inequality holds for sequences $\{x_i(k)\}$ generated by algorithm \eqref{alg:main}:
	\begin{align*}
	\liminf_{k\rightarrow \infty}f(x_{\rm av}(k))\leq f^*+ m\delta
	\end{align*}
\end{theorem}
\pb
To prove this theorem, we consider \eqref{sys:average}.  For any $x^*\in \mathcal{X}^*$, we denote ${\bar x}(k)=x_{\rm av}(k)-x^*$ and obtain that
\begin{align*}
\begin{split}
|{\bar x}(k+1)|^2&=|{\bar x}(k)-\alpha_k \frac{{\bf 1}^\top g_{\delta}(k)}{m}|^2\\
&=|{\bar x}(k)|^2-2\alpha_k\frac{{\bf 1}^\top g_{\delta}(k)}{m} {\bar x}(k)+ |\alpha_k \frac{{\bf 1}^\top g_{\delta}(k)}{m}|^2
\end{split}
\end{align*}
By evaluating \eqref{eq:inexact-distributed} at ${\bf y}={\bf 1} x^*$ and ${\bf y}={\bf 1}x_{\rm av}(k)$ respectively, the following inequalities hold:
\begin{align*}
\begin{split}
{\bf 1}^Tg_{\delta}(k){\bar x}(k)\geq -f(x^*)+(\sum_{i=1}^mf_{i,\delta}(k)-\sum_{i=1}^m \Delta_i(k))\\
f(x_{\rm av}(k))-(\sum_{i=1}^mf_{i,\delta}(k)-\sum_{i=1}^m \Delta_i(k))\leq L\|{\bar {\bf x}}(k)\|^2+m\delta
\end{split}
\end{align*}

Then, under Assumption \ref{ass:boundedness}, we have  
\begin{align*}
\begin{split}
	|{\bar x}(k+1)|^2&\leq |{\bar x}(k)|^2+C^2 \alpha^2_k+\frac{2\alpha_k}{m}[f(x^*)-f(x_{\rm av}(k))+L\|{\bar {\bf x}}(k)\|^2+m\delta] \\
	&\leq |{\bar x}(k)|^2+C^2 \alpha^2_k-\frac{2\alpha_k}{m}[-L\|{\bar {\bf x}}(k)\|^2 +f(x_{\rm av}(k))-m\delta-f(x^*)]
\end{split}
\end{align*}

We claim that $\liminf _{k \to \infty}(f(x_{\rm av}(k))-m\delta-f(x^*))\leq 0$ and try to seek a contradiction. For the purpose, we denote $b_0=\liminf _{k \to \infty}(f(x_{\rm av}(k))-m\delta-f(x^*))$ and assume $b_0>0$. Since $\|{\bar {\bf x}}(k)\|$ goes to zero when $k \to \infty$, there must be a large integer $K>0$ such that $f(x_{\rm av}(k))-m\delta-f(x^*)-L\|{\bar {\bf x}}(k)\|^2>\frac{b_0}{2}$ for any $k>K$. This implies that
\begin{align*}
|{\bar x}(k+1)|^2&\leq |{\bar x}(k)|^2+C^2 \alpha^2_k-\frac{ \alpha_k}{m}b_0
\end{align*}
Or equivalently,
\begin{align*}
\frac{b_0}{m}\sum_{k'=K}^k\alpha_k' \leq |{\bar x}(K)|+C^2\sum_{k'=K}^k\alpha_{k'}^2
\end{align*}
for any $k>K$, which is impossible for a sufficient large $k$ under theorem conditions. Hence, we have found a contradiction.  It is safe for us to conclude that $b_0\leq  0$. The proof is thus complete. 
\pe

\begin{remark}
	 From Theorem \ref{thm:inexact}, we conclude that the distributed optimization problem \eqref{prop:main} is approximately solved by rule \eqref{alg:main} in the sense that $\lim_{k\to\infty}[x_i(k)-x_{\rm av}(k)]=0$ and $\liminf_{k\rightarrow \infty}f(x_{\rm av}(k))\leq f^*+ m\delta$. On the one hand, this quantified relationship is consistent with the centralized result in \cite{nedic2010effect, devolder2014first}. On the other hand, when the first-order oracle is exact with $\delta=0$, we recover the convergence and optimality conclusions obtained in existing literature \cite{nedic2009distributed, jakovetic2014fast, qu2018harnessing}. 
\end{remark}

\begin{remark}
	It is remarkable that the oracle inexactness does not accumulate over time according to the rule \eqref{alg:main}. In fact, it can be verified that only the long-term accuracy matters in the steady-state performance of our algorithm if the oracle inexactness is time-varying. This property might be favorable in real-time noisy circumstances. 
\end{remark}

Next, we will specify some conditions on the oracle accuracy under which the exact convergence of sequences $\{x_i(k)\}$ to an optimal solution is guaranteed. Since diminishing oracle accuracy is required to ensure the exact (consensus) solvability of our problem by Theorem \ref{thm:inexact}, we make another assumption to replace Assumption \ref{ass:inexact}.
\begin{ass}\label{ass:inexact-varying}
	For each $i=1,\,\dots,\,m$, function $f_i$ is continuous and admits a first-order $(\delta_k,\,L_k)$-oracle and $\max_{k}\{L_k\}\leq L$ for some constant $L>0$. 
\end{ass}

Here is the second main theorem of this paper.

\begin{theorem}\label{thm:exact}
Suppose Assumptions \ref{ass:boundedness}--\ref{ass:inexact-varying} hold  and the positive step size satisfies 
	\begin{align*}
		\sum_{k=1}^\infty \alpha_k= \infty, \quad  \sum_{k=1}^\infty \alpha_k^2< \infty, \quad \sum_{k=1}^\infty \alpha_k\delta_k<\infty
	\end{align*}
Then the following assert will take place:
	\begin{itemize}
			\item[a)] The sequence $\{\|x_{\rm av}(k)-x^*\|\}$ converges;
			\item[b)] There exists a cluster point ${\bar x}_{\rm av}$ of sequence  $\{x_{\rm av}(k)\}$ such that ${\bar x}_{\rm av}\in X^*$;
			\item[c)] If the optimal solution to problem \eqref{prop:main} is unique, then $\lim_{k\rightarrow \infty}x_{i}(k)=x^*$ for each $i=1,\,\dots,\,m$. 
	\end{itemize}
\end{theorem}
\pb
Following the proof of Theorem \ref{thm:inexact}, we have 
\begin{align}\label{eq:key-exact}
	\begin{split}
		|{\bar x}(k+1)|^2&\leq |{\bar x}(k)|^2+C^2 \alpha^2_k-\frac{2\alpha_k}{m}[-L\|{\bar {\bf x}}(k)\|^2 +f(x_{\rm av}(k))-m\delta_k-f(x^*)]
	\end{split}
\end{align}
with ${\bar x}(k)=x_{\rm av}(k)-x^*$. Ignoring the negative terms and using the fact $m \geq 1$, we have
\begin{align*}
	\begin{split}
		|{\bar x}(k+1)|^2&\leq |{\bar x}(k)|^2+C^2 \alpha^2_k+ 2 L \alpha_k  \|{\bar {\bf x}}(k)\|^2+2\alpha_k\delta_k
	\end{split}
\end{align*}
We claim that $\sum_{k=1}^{\infty} \alpha_k  \|{\bar {\bf x}}(k)\|^2<\infty$. 

We first show $\sum_{k=1}^{\infty} \alpha_k  \|{\bar {\bf x}}(k)\|<\infty$.  From the proof of Lemma \ref{lem:consensus}, the following inequality holds for any $k\geq 3$: 	
\begin{align*}
\alpha_{k} ||\bar {\bf x}(k)||&\leq m \alpha_{k} q^{k-2} ||\bar {\bf x}(1)||+ \alpha_{k}\alpha_{k-1}\sqrt{m} C||R||+m^{3/2}  C ||R||\alpha_{k} \sum_{k'=1}^{k-2} \alpha_{k'} q^{k-k'-2}
\end{align*}
Completing the square to handle the $\alpha_{k}$-related terms, one can further obtain that 
\begin{align*}
	\alpha_{k} ||\bar {\bf x}(k)||&\leq m^2  ||\bar {\bf x}(1)||^2 \alpha^2_{k} + q^{2(k-2)} +m C||R||( \alpha_{k}^2+ \alpha_{k-1}^2) \\
	&+m^2  C ||R|| \alpha_{k}^2 \sum_{k'=1}^{k-2}q^{k-k'-2}+m^2  C ||R|| \sum_{k'=1}^{k-2} \alpha_{k'}^2 q^{k-k'-2}
\end{align*}
where we use the fact $\sqrt{m}\leq m$ for $m\geq 1$. Summing up all from $k=3$ and rearranging some terms, we have
\begin{align*}
	\sum_{k=3}^\infty \alpha_{k} ||\bar {\bf x}(k)||&\leq ( m^2  ||\bar {\bf x}(1)||^2+2mC||R||+\frac{m^2  C ||R|| }{1-q})\sum_{k=1}^\infty \alpha_{k}^2\\
	&+\frac{1}{1-q^2}+m^2 C ||R|| \sum_{k=1}^\infty  \sum_{k'=1}^{k} \alpha_{k'}^2 q^{k-k'}
\end{align*}

Since $0<q<1$ and $\sum_{k=1}^\infty \alpha_k^2<\infty$, it must hold that $\sum_{k=1}^\infty  \sum_{k'=1}^{k} \alpha_{k'}^2 q^{k-k'}<\infty$ by Lemma 7 in \cite{nedic2010constrained}. Thus, $\sum_{k=1}^{\infty} \alpha_k  \|{\bar {\bf x}}(k)\|<\infty$. Recalling the fact that $\sum_{k=1}^\infty \alpha_k^2< \infty$ implies $ \lim_{k\to \infty} \alpha_k=0$, one can conclude that ${\bar {\bf x}}(k)$ will converge to $0$ by Lemma \ref{lem:consensus}. This jointly with $\sum_{k=1}^{\infty} \alpha_k  \|{\bar {\bf x}}(k)\|<\infty$ means the preceding claim is indeed true. 

At this stage, we can directly use Lemma 2.1 in \cite{Auslender2004InteriorGradient} to conclude the convergence of $|{\bar x}(k)|^2$ (i.e., item a)) and $\sum_{k=1}^\infty \alpha_{k}[f(x_{\rm av}(k))-f^*]<\infty$ under the theorem conditions. 

Since $|x_{\rm av}(k+1)-x_{\rm av}(k)|=|\alpha_k \frac{{\bf 1}^\top g_{\delta}(k)}{m}|\leq C \alpha_k$ under Assumption \ref{ass:inexact-varying}, we resort to Proposition 2 in \cite{alber1998Oprojected} and have $\lim_{k\to \infty} f(x_{\rm av}(k))=f^*$.  At the same time, the sequence $\{x_{\rm av}(k)\}$ is  bounded by item a) and must have a convergent subsequence $\{{\bar x}_{k_{k'}}\}$. We denote the limit of this subsequence by ${\bar x}_{\rm av}$. Then, from the continuity of function $f$, we have $f({\bar x}_{\rm av})=f(\lim_{k'\to\infty}{\bar x}_{k_{k'}})=\lim_{k'\to\infty}f({\bar x}_{k_{k'}})=f^*$, which means ${\bar x}_{\rm av} \in X^*$ (i.e., item b)). 

When the optimal solution is unique, all convergent subsequences of  $\{x_i(k)\}$ converge to the same limit. Then, $\lim_{k\rightarrow \infty}x_{i}(k)=\lim_{k\rightarrow \infty}x_{\rm av}(k)=x^*$ for each $i=1,\,\dots,\,m$. The proof is thus complete. 
\pe

\begin{remark}
	This theorem states that the distributed rule \eqref{alg:main} can generate the exact optimal solution to our problem even with some inexact first-order information provided that the accumulated oracle inexactness is finite in terms of $\sum_{k=1}^\infty \alpha_k\delta_k<\infty$. This condition is motivated by our previous work \cite{zhu2022primal} on distributed primal-dual $\varepsilon$-subgradient method for problem \eqref{alg:main}. Similar conditions have been widely investigated to handle the inexactness of subgradient and proximal point calculation in centralized optimization scenarios \cite{rockafellar1976monotone, nedic2010effect, kiwiel2004convergence, rasch2020inexact}.  This property can be very useful especially when the first-order information is numerically obtained by solving another auxiliary optimization problem with controlled accuracy. 
\end{remark}

\section{Simulation}\label{sec:simu}

\begin{figure}
	\centering
	\subfigure[]{
		\includegraphics[width=0.42\linewidth]{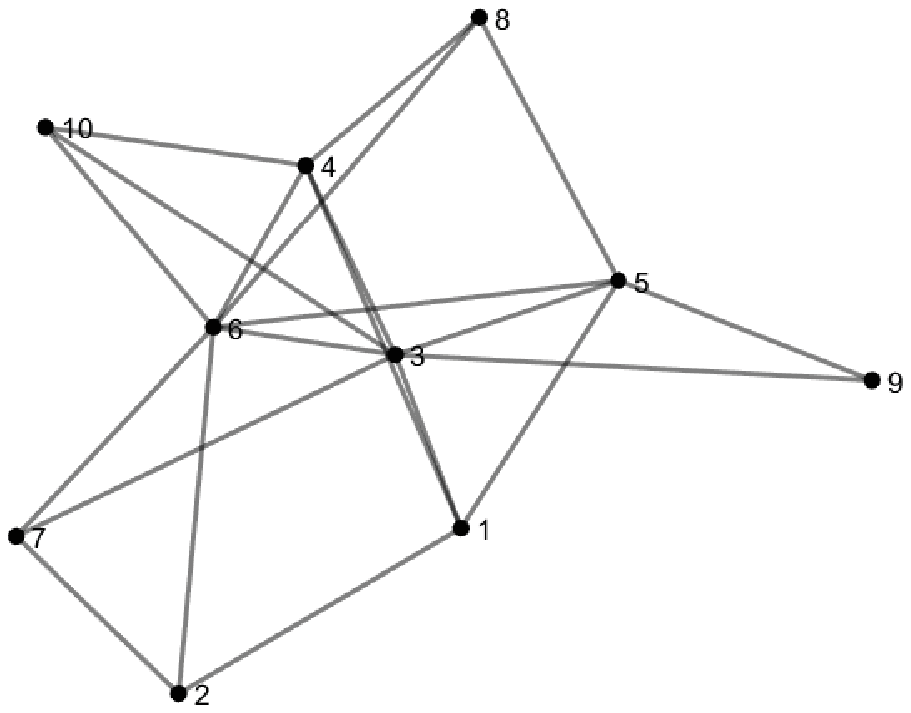}
	} 
	\subfigure[]{
		\includegraphics[width=0.42\linewidth]{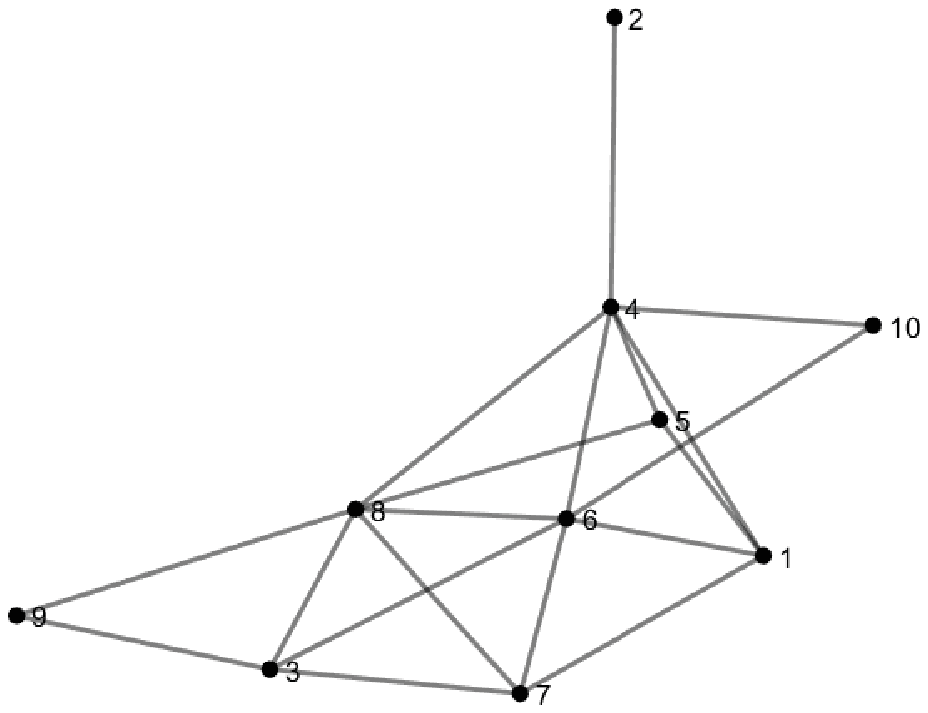}
	}
	\caption{Communication graphs in our example.}\label{fig:graph}
\end{figure}

\begin{figure}
	\centering
	\includegraphics[width=0.84\linewidth]{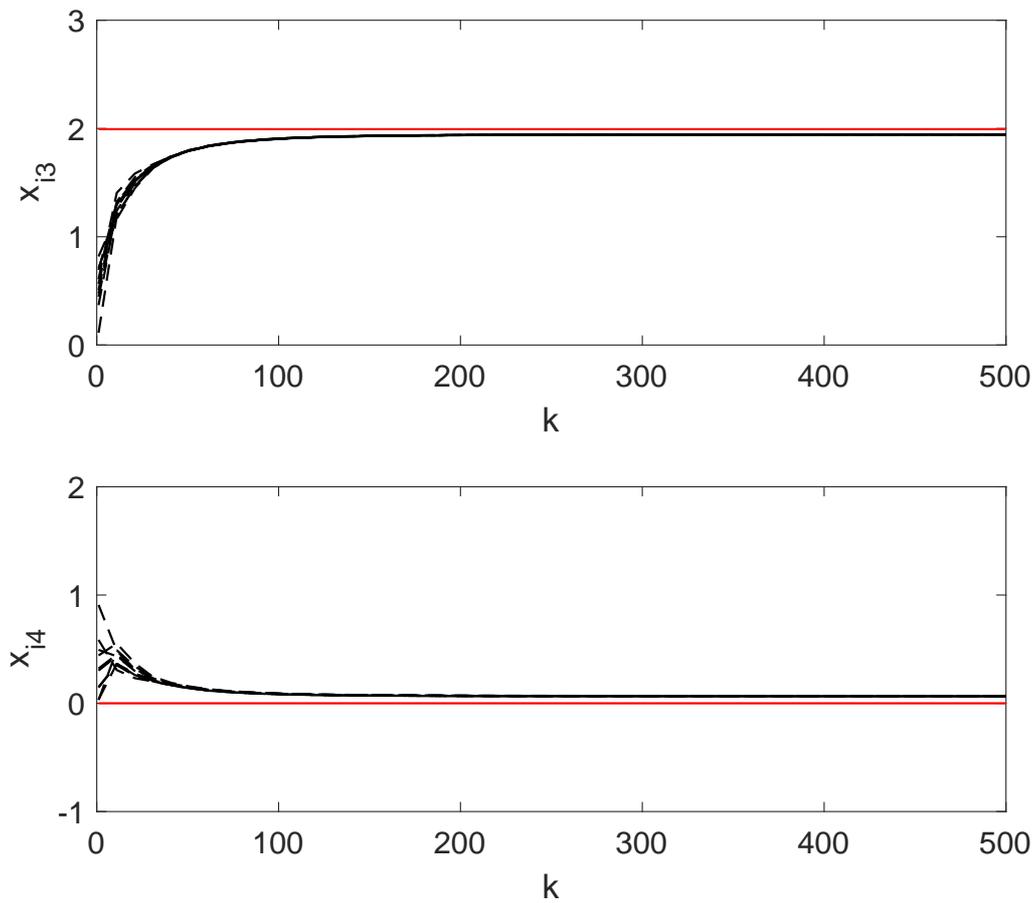}
	\caption{Profile of $x_{i3}$ and $x_{i4}$ with $\delta_k=1$. The estimates are illustrated by black dashed lines and the global optimal solution by red solid  lines. }\label{fig:fixed-error}
\end{figure}

In this section, we provide a numerical example to illustrate the effectiveness of our proposed algorithm in simulations. 

Consider the following LASSO regression problem
\begin{align*}
	\min_{x\in \R^n} \{\frac{1}{2} ||Ax-y||^2+\lambda ||x||_1\}
\end{align*}
where $A\in \R^{N\times n}$, $y\in \R^{N}$,  and $\lambda>0$. Suppose the training data $(y,\, A)$ is collected by $m$ agents in the form:  $A=[A_1,\,\dots,\,A_m]^\top$ and $y=[y_1,\,\dots,\, y_m]^\top$ with $y_i\in \R^{N_i}$, $A_i\in \R^{N_i\times n}$, and $\sum_{i=1}^m N_i=N$. The problem can be rewritten into the form \eqref{prop:main} with $f_i(x)=\frac{1}{2}||A_ix-y_i||^2+\frac{\lambda}{m} ||x||_1$ ($i=1,\,\dots,\,m$) and solved by the distributed subgradient method \cite{nedic2009distributed, liu2017convergence}. Here we use this example to show the effectiveness of our algorithm working with inexact first-order oracle. 

\begin{figure} 
	\centering
	\includegraphics[width=0.84\linewidth]{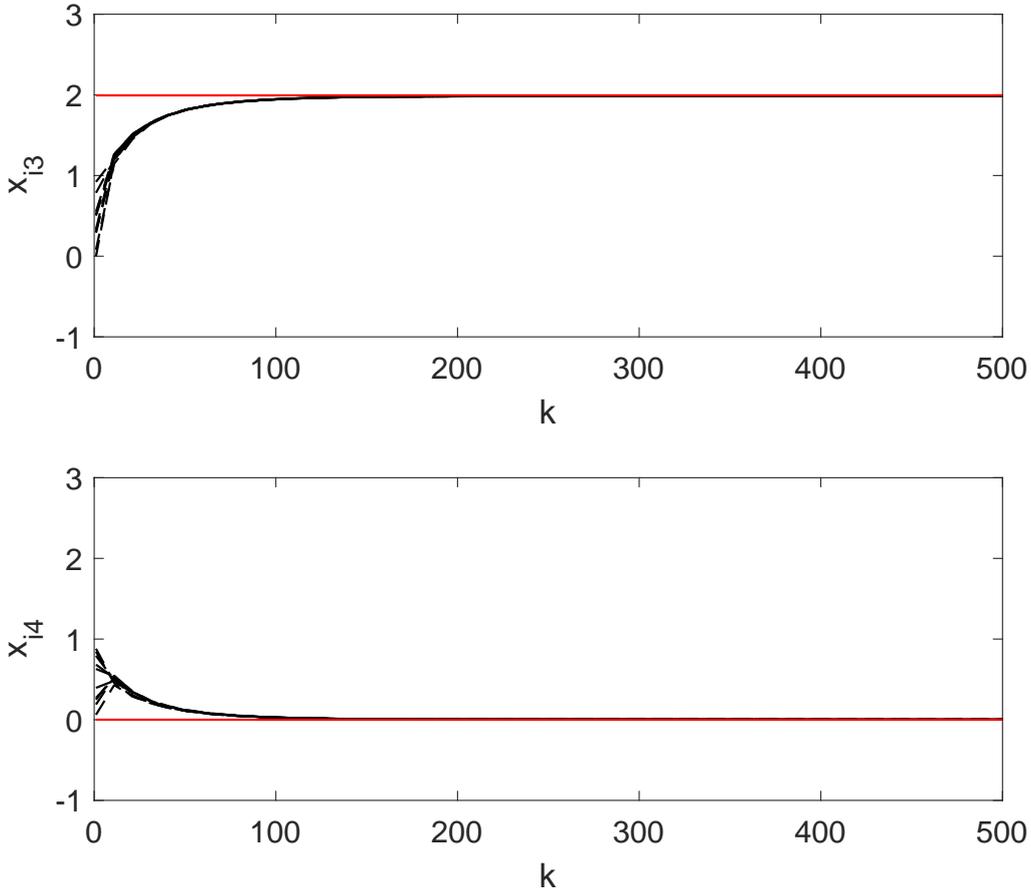}
	\caption{Profile of $x_{i3}$ and $x_{i4}$ with $\delta_k=1/k$. The estimates are illustrated by black dashed lines and the global optimal solution by red solid lines. } \label{fig:varying-error}
\end{figure}

Note that $|x|$ can be lower approximated by the following Huber loss function:
$$ 
H_\delta(x)=
\begin{cases}
	{x^2}/(2 \delta), & |x|\leq \delta\\
	|x|-{\delta}/{2},  & |x|>\delta
\end{cases}   
$$
We denote $\tilde f_{i\delta}(x_i)=\frac{1}{2}||A_ix_i-y_i||^2+\frac{\lambda}{m} \sum_{i=1}^n H_\delta(x_{ij})$ for each $x_i=\mbox{col}(x_{i1},\,\dots,\,x_{in})\in \R^n$ with $0<\delta\leq 1$. It can be verified that $\tilde f_{i\delta}(x_i)$ is convex, differentiable with Lipschitz gradient. Moreover, $(\tilde f_{i\delta}(x_i),\,\nabla \tilde f_{i\delta}(x_i))$ is an $(n\delta/2,\, \sqrt{n}/\delta+||A_i||^2)$-oracle for the nonsmooth objective function $f_i(x_i)$.  Then we can apply the distributed rule \eqref{alg:main} to resolving this LASSO regression problem.

\begin{figure}
	\centering
	\includegraphics[width=0.84\linewidth]{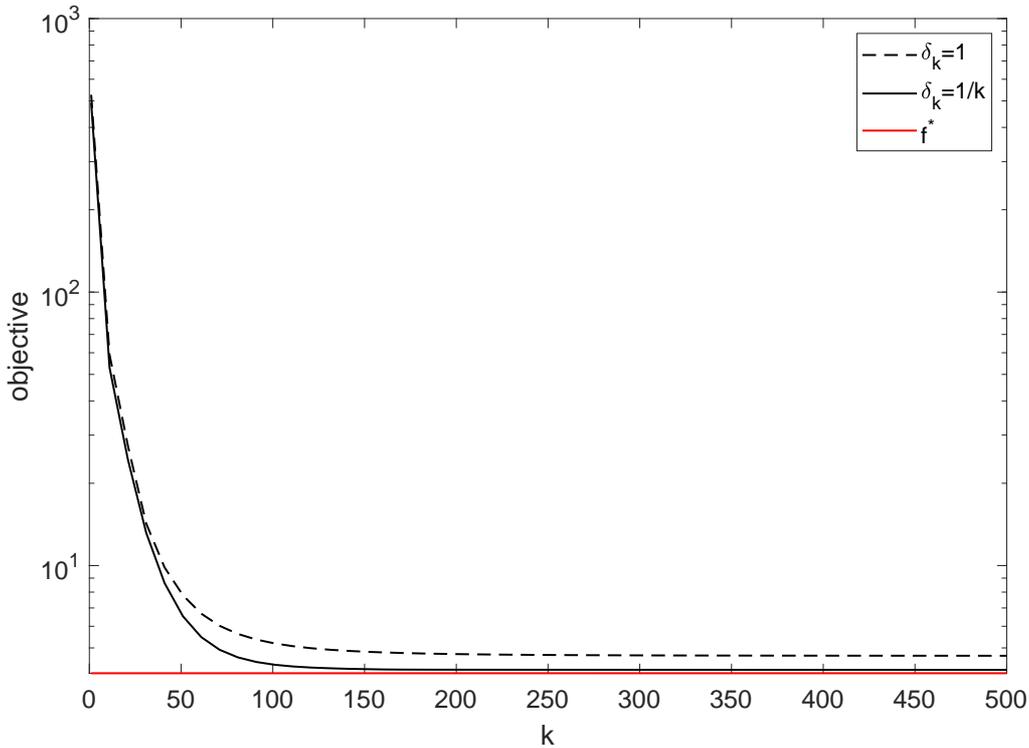}
	\caption{Profile of $\sum_{i=1}^m \tilde f_{i\delta}(x_i)$ with different inexact oracles.}\label{fig:f-compare}
\end{figure}

In simulations, we let $N=1000$, $m=10$, $N_i=100$, $n=5$ and $\lambda=1$.  The matrix $A$ is randomly chosen with each entry in $[0,\,1]$ and $y=Ax_{0}+\epsilon$ with $x_{0}=\mbox{col}(0,\, 1,\, 2,\, 0,\, 1)$ and $\epsilon$ randomly chosen in $[-0.01,\,0.01]$.  Suppose the communication topology alternates between the two graphs depicted in Figure \ref{fig:graph} every 50 iterations. The corresponding two weighted matrices are randomly chosen to satisfy Assumption \ref{ass:graph-weight}.
The step size is chosen as $\alpha_{k}=0.5/(k+80)$ to ensure a consensus of all agents' estimates.  For the inexact oracle, we first fix the oracle accuracy at $\delta_k=1$ and then use the time-varying one with $\delta_k=1/k$. The profiles of two components of agent $i$'s estimate $x_i$ (i.e., $x_{i3}$ and $x_{i4}$) under different oracle conditions are reported in Figures \ref{fig:fixed-error} and \ref{fig:varying-error}. In Figure \ref{fig:fixed-error}, with constant oracle accuracy,  the estimates are observed to reach a consensus but deviate from the global optimal solution $x^*$.  By contrast, these estimates converge to $x^*$ when we use inexact oracle with time-varying accuracy $\delta_k$ as depicted in Figure \ref{fig:varying-error}.  The profile of $ \sum_{i=1}^m \tilde f_{i\delta}(x_i)$ in two cases and also the optimal value are shown in Figure \ref{fig:f-compare}. It can be found that the distributed rule \eqref{alg:main} can indeed produce an approximation of the optimal value when $\delta_k=1$ while the term $\sum_{i=1}^m \tilde f_{i\delta}(x_i)$ eventually converges to $f^*$ when $\delta_k=\frac{1}{k}$. These observations confirm the derived theoretical conclusions on the proposed rule \eqref{alg:main} to solve the distributed problem \eqref{prop:main} with inexact first-order information.

\section{Conclusion}\label{sec:con}
This paper has attempted to solve a distributed optimization problem using only the approximate first-order information of the objective functions. We have proposed a variant of the distributed subgradient method working with inexact oracle and discussed its convergence properties under different stepsize and oracle conditions. In the future, we may consider the convergence rate of the iterative sequences and further exploit its deep connection with other classes of distributed optimization algorithms.

\bibliographystyle{IEEEtran}
\bibliography{inexact_j}

\end{document}